\documentclass[letterpaper, 10 pt, conference]{ieeeconf}  

\IEEEoverridecommandlockouts                              

\usepackage{graphicx} 
\usepackage{amsmath} 
\usepackage{amssymb}  
\usepackage{color}
\renewcommand{\baselinestretch}{0.99}
\DeclareMathOperator{\e}{\mathrm{e}}

\title{\LARGE \bf
Nonlinear dynamics in pulse-modulated feedback drug dosing\textsuperscript{*}
}

\author{Alexander Medvedev, Anton V. Proskurnikov, and Zhanybai T. Zhusubaliyev
\thanks{*AM was partially supported the Swedish Research Council under grants 2019-04451, 2024-04943.}
\thanks{Alexander Medvedev [{\tt\small alexander.medvedev@it.uu.se}] is with Department of Information Technology, 
        Uppsala University, SE-752 37 Uppsala, Sweden.
        }%
\thanks{Anton V. Proskurnikov  [{\tt\small anton.p.1982@ieee.org}] is with Department of Electronics and Telecommunications, Politecnico di Torino, Turin, Italy, 10129.}%
\thanks{Zhanybai T. Zhusubaliyev [{\tt\small zhanybai@hotmail.com}] is with Department of Computer Science, International Scientific Laboratory for
Dynamics of Non-Smooth Systems, Southwest State University, Kursk, Russia and Faculty of Mathematics and Information Technology, Osh State University, Lenin st. 331, 723500, Osh, Kyrgyzstan.}
}

\begin{document}

\maketitle
\thispagestyle{empty}
\pagestyle{empty}

\begin{abstract}
Pulse-modulated feedback is utilized in drug dosing to mimic sustained over a longer period of time manual discrete dose administration, the latter is in contrast with continuous drug infusion. The intermittent mode of dosing calls for a hybrid (continuous-discrete) modeling of the closed-loop system, where the pharmacokinetics and pharmacodynamics of the drug are captured by differential equations whereas the control law is described by difference equations. Hybrid dynamics are highly nonlinear which complicates formal design of pulse-modulated feedback. This paper demonstrates complex nonlinear dynamical phenomena arising in a simple control system of dosing a neuromuscular blockade agent in anesthesia. Along with the nominal periodic regimen, undesirable nonlinear behaviors, i.e. periodic solutions of high multiplicity, multistability, as well as deterministic chaos, are shown to exist. It is concluded that design of feedback drug dosing algorithms based on a hybrid paradigm has to be informed by a thorough bifurcation analysis in order to secure patient safety.

{\textbf{\textit{Clinical Relevance}}}\textemdash The paper highlights potential patient safety risks posed by complex nonlinear phenomena in closed-loop drug administration. It also presents a systematic  approach to mitigate them at the controller design phase.

\end{abstract}

\section{INTRODUCTION}
 Relieving the burden of manual drug administration over lengthy periods of time requires automation. Surgery with  operative time over two hours is common and has to be reliably supported by general anesthesia that allows patients to be unconscious and free of pain throughout the procedure. 
 
 General anesthesia is nowadays predominantly  achieved by intravenous (IV) administration of
sedatives (hypnotics), analgesics, and muscle relaxants, i.e. neuromuscular blockade (NMB) agents.  Anesthetic drugs are given intravenously  either as  bolus (one-time) doses  or a continuous infusion. A standard procedure in surgery  consists of an initial NMB bolus followed by continuous infusion of the drug at a constant flow rate. This mode of administration is supported in target-controlled infusion (TCI), a computerized system that determines the infusion rate required to produce a desired drug concentration at the effect site. In TCI,  the required dose is calculated using a pharmacokinetic-pharmacodynamics (PKPD) drug model from patient's age, sex, weight and other relevant parameters.

Programmed intermittent bolus (PIB) is an anesthesia technique in which boluses of an anesthetic drug are automatically injected multiple times, with or without patient-controlled boluses.  Being delivered into the epidural space in labor analgesia, PIB has shown to reduce local anesthetic usage and improve maternal satisfaction \cite{GAH13}. The lack of reliable real-time pain sensing technology \cite{FBW23} hinders  development of a closed-loop version of PIB, where bolus delivery would be controlled by an objective effect measurement.

In most cases, NMB agents are administered via continuous infusion.
However,  controlled boluses are recommended in some cases because intermittent doses allow serial evaluation and
reduce the risk of developing myopathies due to prolonged paralysis \cite{RRR22}. The effect of NMB agents is routinely measured by neuromuscular monitors~\cite{MH06}, devices that electrically stimulate a peripheral nerve while also quantifying the evoked responses. Compared to the administration of fixed doses (open-loop control), using the monitors for dose titration during the course of treatment significantly reduces the exposure to NMB drugs without affecting the observed clinical outcome~\cite{TSB21}.

The present paper investigates the dynamical properties of a drug dosing system that implements PIB administration in a feedback framework. The NMB is selected as application due to the availability of reliable effect quantification that enables feedback control.
Using a mixture of analytical and numerical analysis methods, it is demonstrated that complex nonlinear dynamical phenomena can arise in the closed-loop drug administration when continuous effect measurement is utilized to control the regimen of intermittent boluses.

The paper is organized as follows. First mathematical models adopted in the study are introduced. Then, the dynamical properties of the closed-loop drug delivery system are summarized. Based on a parsimonious PKPD model of a NMB agent, a pulse-modulated feedback law implementing a clinically established dosing regimen is designed. Further, bifurcation analysis of the closed-loop dynamics is performed to identify complex nonlinear phenomena arising in the drug dosing system operation and illustrate them by simulation. Finally, conclusions are drawn.

\section{Mathematical  models}\label{sec:nonlin}
\paragraph*{Continuous part}
Consider a time-invariant Wiener system whose measured output is a nonlinear function of the linear block output. The linear block is given by the state-space representation
\begin{equation}                            \label{eq:1_wiener}
\dot{x}(t) =Ax(t), \quad \bar y(t)=Cx(t),
\end{equation}
where the matrices are
\begin{equation}                            \label{eq:1a}
A=\begin{bmatrix} -a_1 &0 &0 \\ g_1 & -a_2 &0 \\ 0 &g_2 &-a_3 \end{bmatrix}, 
C =\begin{bmatrix}0 & 0 & 1\end{bmatrix},
\end{equation}
 $a_1,a_2,a_3>0$ are distinct constants, and $g_1,g_2>0$ are positive gains. The measured output is then
\begin{equation}\label{eq:nonlin}
y(t)=\varphi(\bar y),
\end{equation}
where $\varphi(\cdot)$ is a smooth function.
\paragraph*{Discrete part}
Continuous-time system~\eqref{eq:1_wiener} is controlled by a pulse-modulated feedback that gives rise to instantaneous  jumps in the state vector $x(t)$ 
\begin{align}\label{eq:2}
x(t_n^+) &= x(t_n^-) +\lambda_n B, \quad                                  
t_{n+1} =t_n+T_n, 
\\ 
 B&=\begin{bmatrix} 1 & 0 & 0\end{bmatrix}^{\top},\quad n=0,1,\ldots, \notag
\end{align}
where 
\[
T_n =\bar\Phi(y(t_n)), \; \lambda_n= \bar F(y(t_n)).
\]
In pulse-modulated control, $\bar\Phi(\cdot)$ is  referred to as the frequency modulation function and $\bar F(\cdot)$  as the amplitude modulation function. 
The minus and plus in a superscript in~\eqref{eq:2} denote the left-sided and
right-sided limits, respectively. The described by \eqref{eq:2} control mechanism corresponds to plant \eqref{eq:1_wiener} being subject to an impulsive action $\lambda_n B\delta(t_n)$ applied directly to the state vector, where $\delta(\cdot)$ is Dirac delta function. Despite the jumps in \eqref{eq:2}, the output $y(t)$ remains a smooth function  since $\varphi(\cdot)$ is continuous and
 \begin{equation}\label{CBLB}
 CB=CAB=0,\,CA^2B\ne 0.
 \end{equation}
Since both $\lambda_n$ and $T_n$ are functions of the output value at $t_n$, control law \eqref{eq:2} implements a discrete first-order feedback loop over third-order continuous plant \eqref{eq:1_wiener}.

With $\circ$ denoting composition, introduce the functions
\begin{equation*}
\Phi(\cdot)\triangleq (\bar\Phi\circ\varphi)(\cdot), \quad F(z)\triangleq (\bar F\circ\varphi)(\cdot).
\end{equation*}
Then closed-loop system \eqref{eq:1_wiener}, \eqref{eq:2} is equivalent to the Impulsive Goodwin's Oscillator model \cite{MCS06},\cite{Aut09}. 
It constitutes a hybrid (continuous-discrete) system that is able to exhibit a wide range of nonlinear dynamics phenomena but can also be designed to produce a desired behavior through the choice of the modulation functions $\bar F(\cdot)$, $\bar \Phi(\cdot)$. 

\subsection{Pharmacokinetic-pharmacodynamic model}\label{sec:PKPD}
A minimally parametrized pharmacokinetic-pharmacody\-namic (PKPD) model of the neuromuscular blockade  (NMB) agent \emph{atracurium} is introduced in \cite{SWM12}. It captures the dose-effect relationship for the drug with a linear dynamical block representing the pharmacokinetics and a nonlinear (static) function describing the pharmacodynamics.
Being rewritten in state space, it has the form of \eqref{eq:1_wiener}, \eqref{eq:nonlin} with a certain choice of the parameters and the  output function $\varphi(\cdot)$. 

The output $y(t)$ $\lbrack \% \rbrack$ represents  the effect of the NMB agent and is measured by a train-of-four neuromuscular monitor~\cite{MH06}. The maximal level of $y(t)=100\%$ is achieved when the  NMB is initiated and there is no drug in the bloodstream of the patient. 
 The elements of the matrix $A$, cf. \eqref{eq:1a}, are  parametrized in terms of a common factor $\alpha>0$ 
\begin{equation}\label{eq:A_NMB}
    a_1\!=\! v_1\alpha, a_2\!=\! v_2\alpha, a_3\!=\! v_3\alpha, g_1\!=\! v_1\alpha, g_2\!=\! v_2v_3\alpha^2,
\end{equation}
where $v_1=1$, $v_2=4$, $v_3=10$ are fixed coefficients calculated from clinical data. The PD part is modeled by a Hill function of order $\gamma$
\begin{equation}\label{eq:nonlin_NMB}
y=\varphi(\bar y)=\frac{100 C_{50}^\gamma}{ C_{50}^\gamma + {\bar y}^\gamma(t)}, \quad \gamma>0.
\end{equation}
Here $C_{50}=3.2425$ $ \mu \mathrm{ g} \ \mathrm{ml}^{-1} $ is the drug concentration that produces 50\% of the maximum effect.

The minimality of  the PKPD model consists  in its particularization to a patient by means of the pair of constants $(\alpha, \gamma)$.
When the model at hand is estimated from clinical data from 48 patients  \cite{SWM12}, the following parameter intervals are obtained
$ \alpha_{\min}= 0.0270\le \alpha \le 0.0524=\alpha_{\max}$, $\gamma_{\min}=1.4030\le\gamma\le 5.5619=\gamma_{\max}$. The population mean values are evaluated to $\bar \alpha=0.0374,
\bar \gamma=2.6677$.

For the population mean values of the patient-specific parameters, the PK (linear continuous) part of the model has the state matrix 
\[
A=\begin{bmatrix}
     -0.0374         &0         &0\\
    0.0374   &-0.1496         &0\\
         0    &0.0560   &-0.3740
\end{bmatrix}.
\]

In terms of control law \eqref{eq:2}, a dose of $\lambda_n~\mathrm{\mu g/kg}$ of {\it atracurium} is delivered at the time instant $t_n$. In surgical practice, after an initial bolus  dose of $400\text{--}500~\mathrm{\mu g/kg}$, maintenance doses of $80\text{--}200~\mathrm{\mu g/kg}$ are  administered every $10\text{--}20~\mathrm{min}$.

It should be emphasized here that the PKPD model in \eqref{eq:1_wiener},\eqref{eq:A_NMB},\eqref{eq:nonlin_NMB} has been developed explicitly for closed-loop anesthesia and complies with the model parsimonicity principle applied in automatic control. Minimally parametrized models allow to keep the variance of parameter estimates law even though the excitation of the plant model is limited in e.g. maintenance phase. Conventional PKPD models are based on population data, involve multiple covariates, and are produced for optimizing design of clinical studies. In control law \eqref{eq:2}, the initial bolus $\lambda_0$ administered at time $t_0$ is not defined by the measured output and can be calculated from patient data using a conventional PKPD model and an established effect target.


\section{Dynamics}
Denote the plant state vector value prior to a feedback firing as $X_n=x(t_n^-)$. As shown in \cite{MCS06}, \cite{Aut09}, the closed-loop dynamics of hybrid system \eqref{eq:1_wiener},\eqref{eq:2} obey the discrete map
\begin{align}\label{eq:map}
    X_{n+1}&=Q(X_n),\\
    Q(\xi) &\triangleq \mathrm{e}^{A\Phi(C\xi)}\left( \xi+ F(C\xi)B \right).\nonumber
\end{align}
 Given the vector $X_n$, the continuous state trajectory in between the feedback firing instants  is calculated as
\begin{equation} \label{eq:1d}
x(t)=\e^{(t-t_n)A}(X_n+\lambda_n B),\quad t\in(t_n,t_{n+1}).
\end{equation}
Therefore, the sequence $X_n, n=0,1,\dots$ completely defines the hybrid  system dynamics. 
The dynamical behaviors of system \eqref{eq:1_wiener},\eqref{eq:2} have been studied in \cite{ZCM12b} by bifurcation analysis for smooth modulation functions under the conditions below.
\begin{description}
    \item[\bf C1:] $F(\cdot)$ and $\Phi(\cdot)$ are continuous and monotonic;
    \item[\bf C2:] $F(\cdot)$ is non-increasing, and $\Phi(\cdot)$ is non-decreasing on $[0,\infty)$;
    \item[\bf C3:] For some constants $\Phi_1$, $\Phi_2$, $F_1$, $F_2$, it holds that 
    \begin{equation}                             \label{eq:mod_bounds}
0<\Phi_1\le \Phi(\cdot)\le\Phi_2, \ 0<F_1\le F(\cdot)\le F_2.
\end{equation}
\end{description}

All the solutions of \eqref{eq:1_wiener},\eqref{eq:2} are positive because the matrix $A$ is Metzler and the amplitude modulation function is positive. The boundedness of the modulation functions ({\bf C3}) ensures that all the closed-loop solutions are bounded, see \cite{Aut09}.  Since both $\Phi_1$ and $F_1$ are required to be positive ({\bf C3}), the system does not possess equilibria and has only oscillating solutions. 

The simplest type of periodic solution is characterized by just one firing of the feedback on the least period and is called 1-cycle. For a 1-cycle, it applies that the map $Q(\cdot)$ has a fixed point $X$ solving the equation
\begin{equation}\label{eq:1-cycle}
    X=Q(X).
\end{equation}
Let the 1-cycle corresponding to the fixed point have the least period of $T$ and the impulse weight $\lambda$. Then, as shown in \cite{PRM24}, there is a closed-form expression  for the fixed point
   \begin{equation}\label{eq:fp_alpha}
        X= \lambda {(\e^{-TA}-I)}^{-1} B,
    \end{equation}  
    and the fixed point is positive.
 It follows from \eqref{eq:fp_alpha}  that 1-cycle always exists and unique, under the assumptions made. Yet, it is not necessarily orbitally stable and can evolve to another solution under an infinitesimal perturbation.

Denote the output of the linear block at the fixed point  $\bar y_0\triangleq CX$. The Jacobian of the pointwise map $Q(\cdot)$ evaluated at the fixed point $X$ is given by
    \begin{equation}\label{eq:jacobian}
    Q^\prime(X)=A_\Phi+WKC,
\end{equation}
\[
A_\Phi=\e^{A\Phi(\bar y_0)}, W=\begin{bmatrix}
    J &D
\end{bmatrix}, K^\top= \begin{bmatrix}
  F^\prime(\bar y_0)   & \Phi^\prime(\bar y_0)
\end{bmatrix}.
\]
A necessary and sufficient condition for a 1-cycle being orbitally stable is that $Q^\prime(X)$ is Schur-stable
\begin{equation}\label{eq:Schur}
    \rho\left( Q^\prime(X) \right)<1,
\end{equation}
where $\rho(\cdot)$ is spectral radius.

In \cite{PRM24}, it is proved that
\[
D<0\;\;\text{and}\;\; J>0,
\]
which property, together with {\bf C2}, gives 
\begin{equation}\label{eq:negative_feedback}
   WK= JF^\prime(\bar y_0) + D\Phi^\prime(\bar y_0) <0.
\end{equation}
Inequality \eqref{eq:negative_feedback} implies that feedback law \eqref{eq:2} implements a negative feedback over plant \eqref{eq:1_wiener} which is well-defined for all feasible values of $F^\prime(\bar y_0),\Phi^\prime(\bar y_0)$, despite the hybrid nature of the underlying closed-loop system.

Periodic solutions of higher multiplicity can be obtained by applying iterations of the map $Q(\cdot)$
\begin{equation}\label{eq:Q_m}
    Q^{(m)}(X)= \underbrace{Q(Q(\dots(Q}_{m}(X))\dots)).
\end{equation}
Then, by the chain rule, the Jacobian of map \eqref{eq:Q_m} is
\begin{align}\label{eq:Q_m_prime}
     {\left(Q^{(m)}(X)\right)}^\prime&=Q^\prime\left(Q^{(m-1)}(X)\right)Q^\prime\left(Q^{(m-2)}(X)\right)\dots\notag \\
     &\times Q^\prime\left(Q(X)\right) Q^\prime(X).
\end{align}

Existence of an $m$-cycle is equivalent to the algebraic system
\begin{equation}\label{eq:m-cycle}
    X=Q^{(m)}(X),
\end{equation}
having a (positive) solution $X$. The $m$-cycle is stable when 
$$\rho\Big({\left(Q^{(m)}(X)\right)}^\prime\Big)<1.$$

Besides orbital stability, the Jacobian matrix communicates other relevant solution properties. The spectral radius of the Jacobian defines the rate at which the solution tends to the  stationary one after a small perturbation from it. Naturally, this convergence rates applies only locally and is not valid for large transients. The eigenvalues of the Jacobian are also the multipliers of the fixed point and (locally) define the character of the transients. To achieve e.g. a transient without overshoot, the modulation functions $F(\cdot)$ and $\Phi(\cdot)$ have to be selected so that the multipliers are positive.

Since the measured output of \eqref{eq:1_wiener} varies significantly, e.g. $y(t)\in \lbrack 0, 100 \rbrack$ in \eqref{eq:nonlin_NMB}, it is essential to the closed-loop system function that all feasible initial conditions $x(t_0), \lambda_0$ result in a solution that belongs to the basin of attraction of the desired stationary solution.

\section{Design}\label{sec:design}

This section briefly goes through the steps of a design algorithm that yields the parameters of the modulation functions for a given by the constants $(\alpha,\gamma)$ PKPD model  and the parameters of a desired 1-cycle. The design procedure is covered in more detail in  \cite{MPZ24_ECC}. 

Consider the pharmacokinetic-pharmodynamic  model in Section~\ref{sec:PKPD} evaluated for the population mean values $\bar\alpha, \bar\gamma$. Select the parametrization of the modulation functions  of controller \eqref{eq:2}  
as piecewise affine, i.e.
\begin{align}\label{eq:affine_Phi}
    \Phi (\xi)= \begin{cases} \Phi_2 &   \Phi_2 < k_2\xi +k_1, \\
     k_2\xi +k_1 & \Phi_1 \le  k_2\xi +k_1 \le \Phi_2, \\
    \Phi_1  &  k_2\xi +k_1 < \Phi_1, 
     \end{cases}
\end{align}
\begin{align}\label{eq:affine_F}
    F (\xi)= \begin{cases} F_1 &  k_4\xi +k_3< F_1, \\
     k_4\xi +k_3 & F_1 \le k_4\xi +k_3 \le F_2, \\
    F_2 & F_2 <k_4\xi +k_3.
     \end{cases}
\end{align}
Then the parameter set $k_1,k_2,k_3, k_4$ and $F_1,F_2$, $\Phi_1,\Phi_2$ completely describes pulse-modulated controller \eqref{eq:affine_F}, \eqref{eq:affine_Phi}. The  limits of the modulation functions, as well as the parameters of the desired (stationary) periodic solution, are derived from clinical practice. The parameters defining transient performance, i.e. $k_1, k_3$, are tuned manually. Notice that  modulation functions \eqref{eq:affine_Phi}, \eqref{eq:affine_F} are not continuous, i.e. do not comply with {\bf C1}, which property leads to more complex closed-loop dynamics than for a continuous case.

The design aim is to render  an orbitally stable  1-cycle with the parameters $\lambda=300\mathrm{\mu g/kg}$, $T=20~\mathrm{min}$ in the closed-loop system comprising \eqref{eq:1_wiener}, \eqref{eq:nonlin_NMB}, \eqref{eq:2}. Using  \eqref{eq:fp_alpha} the fixed point corresponding to the desired periodic solution is 
\begin{equation}\label{eq:fp_L300_T20}
    X^\top=\begin{bmatrix}  269.5974 &84.5819 & 13.6249\end{bmatrix},
\end{equation}
thus yielding $\bar y_0=13.6249$.

With 
\begin{equation}\label{eq:F_Phi_fixed_point}
    F^\prime(\bar y_0)=-0.15,\quad \Phi^\prime(\bar y_0)=0.29,
\end{equation}
 the eigenvalues of the Jacobian in \eqref{eq:jacobian} are $\sigma(Q^\prime(X))=\{0.2288, 0.1863, 0.0003\}$, which guarantees local stability of the fixed point.

Assuming that $F^\prime(\cdot)$ and $\Phi^\prime(\cdot)$ are outside of saturation, the chain rule gives
\begin{align}\label{eq:F_prime_Phi_prime}
     F^\prime(\bar y_0)&=\bar F^\prime(\bar y_0)\varphi^\prime(\bar y_0)= k_4 \varphi^\prime(\bar y_0),\\ 
     \Phi^\prime(\bar y_0)&= \bar \Phi^\prime(\bar y_0)\varphi^\prime(\bar y_0)= k_2 \varphi^\prime(\bar y_0), \nonumber
\end{align}
where
\[
\varphi^\prime(\xi)= -\frac{\gamma 100 C_{50}^\gamma  \xi^{\gamma-1}}{ {(C_{50}^\gamma + {\xi}^\gamma)}^2}, \quad \varphi^\prime(\bar y_0)=-0.4073.
\]
Thus, from \eqref{eq:F_Phi_fixed_point}  and \eqref{eq:F_prime_Phi_prime}, it follows that
\[ k_2=-0.7119 , \quad k_4=0.3682.
\]
Now, the rest of the coefficients of the modulation functions are obtained from the  parameters of the desired 1-cycle
\begin{align*}
    F(\bar y_{0})&= (\bar F \circ \varphi)(\bar y_{0})=\bar F ( \varphi(\bar y_{0}))=k_4 \varphi(\bar y_{0})+k_3= \lambda,\\
    \Phi(\bar y_{0})&= (\bar \Phi \circ \varphi)(\bar y_{0})= \bar \Phi (\varphi(\bar y_{0}))= k_2 \varphi(\bar y_{0})+k_1=T,
\end{align*}
which, for $\varphi(\bar y_0)=2.1256$, yield 
\begin{gather*}
 k_1= 21.5133, \;k_3 = 299.2173.
\end{gather*}
 The  limits of the  modulation functions are selected as
$F_1=150~\mathrm{\mu g/kg}$,  $F_2=400~\mathrm{\mu g/kg}$, $\Phi_1=11~\mathrm{min}$, $\Phi_2=50~\mathrm{min}$.
Therefore, control law \eqref{eq:2} cannot produce a dose higher than $F_2$ or lower than $F_1$. Further, there is at least one dose administered in $\Phi_2$ minutes but not closer than $\Phi_1$ minutes to the previous one.

\begin{figure}[t]
\centering 
\includegraphics[width=0.9\linewidth]{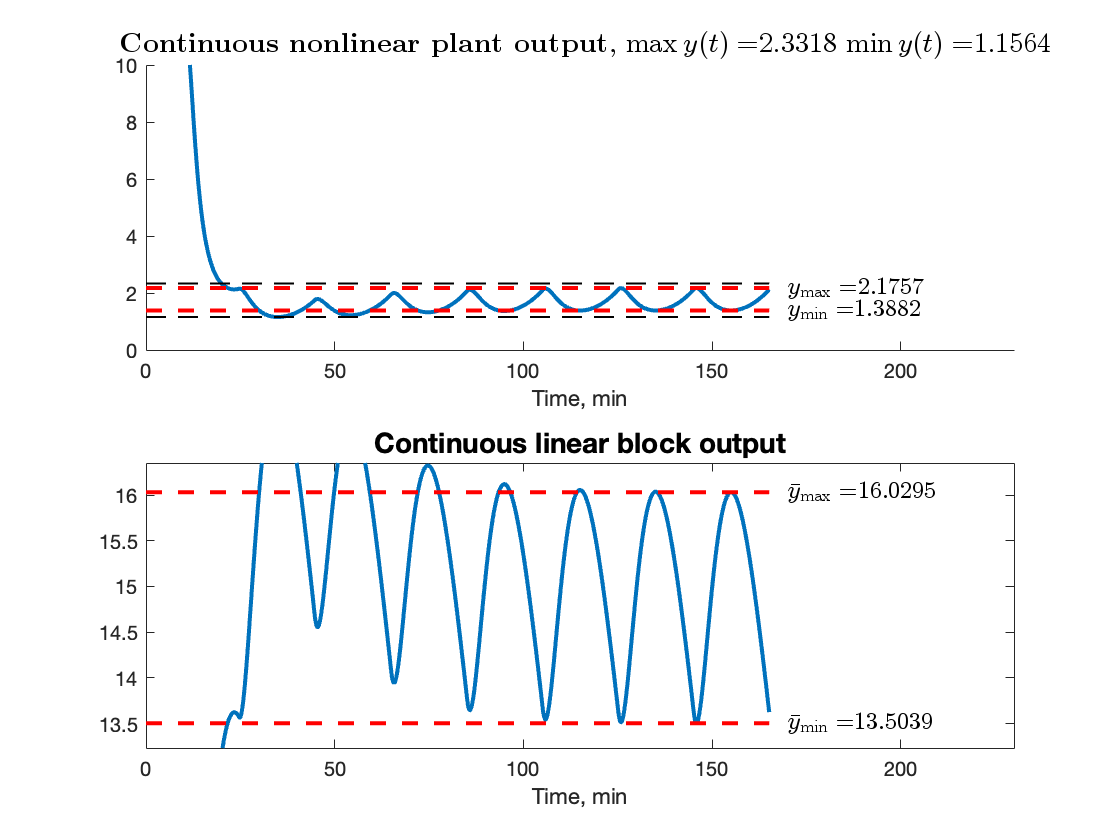}
\caption{ Convergence to the  1-cycle from $x(0)=0$  in the NMB model with $(\bar\alpha,\bar\gamma)$ stabilized by the modulation function slopes $F^\prime(y_0)=-0.15$, $\Phi^\prime(y_0)=0.29$. Top plot: the nonlinear output $y(t)$. The horizontal black dashed lines mark $\inf_t y(t)$ and $\sup_{t\in \lbrack T,5T\rbrack} y(t)$. The stationary output corridor values for  the 1-cycle are marked in red. Bottom plot: the linear output $\bar y(t)$. 
}\label{fig:NMB_1-cycle}
\end{figure}

\begin{figure}[t]
\centering 
\includegraphics[width=0.9\linewidth]{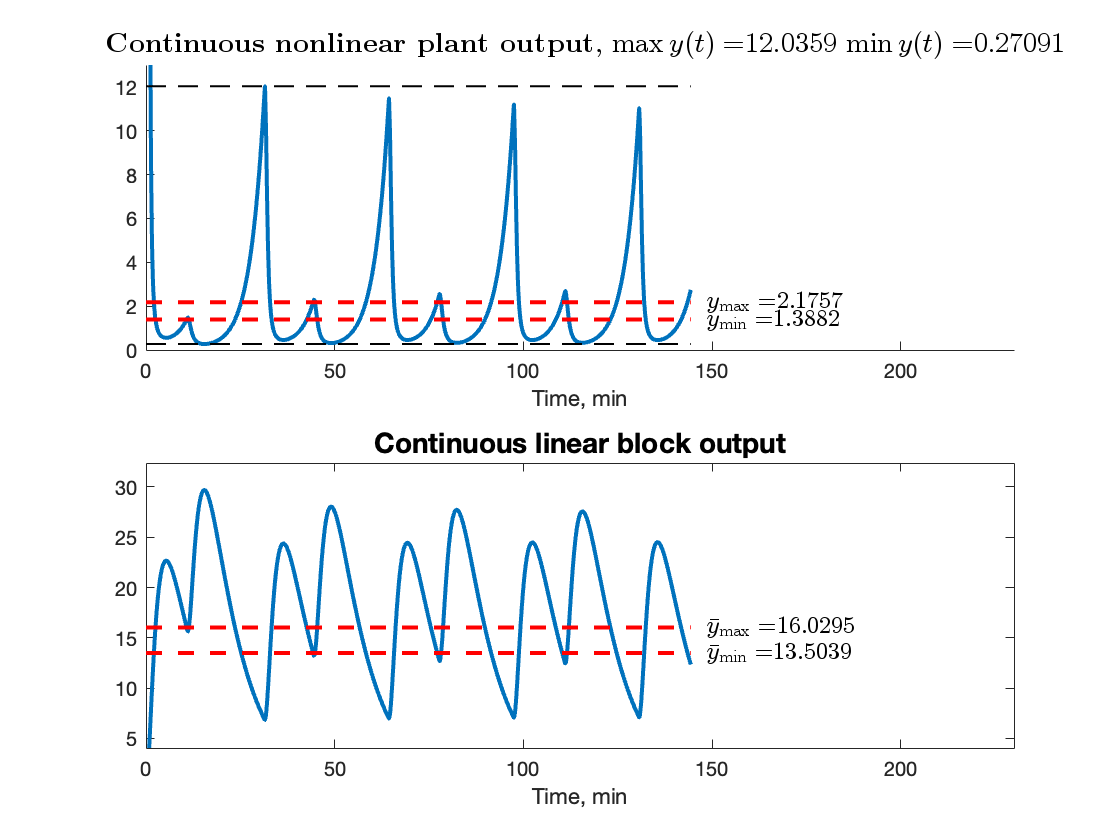}
\caption{ Convergence to the 2-cycle from $x(0)=0$  in the NMB model with $(\alpha=0.11, \bar\gamma)$ stabilized by the modulation function slopes $F^\prime(y_0)=-0.15$, $\Phi^\prime(y_0)=0.29$. Top plot: the nonlinear output $y(t)$. Bottom plot: the linear output $\bar y(t)$. The desired output corridor values for  the 1-cycle are marked in red. The horizontal black dashed lines mark $\inf_t y(t)$ and $\sup_{t\in \lbrack T,5T\rbrack} y(t)$.
}\label{fig:NMB_2-cycle}
\end{figure}

\begin{figure}[ht]
\centering \centering
\parbox[c]{0.65\linewidth}{\includegraphics[width=\linewidth]{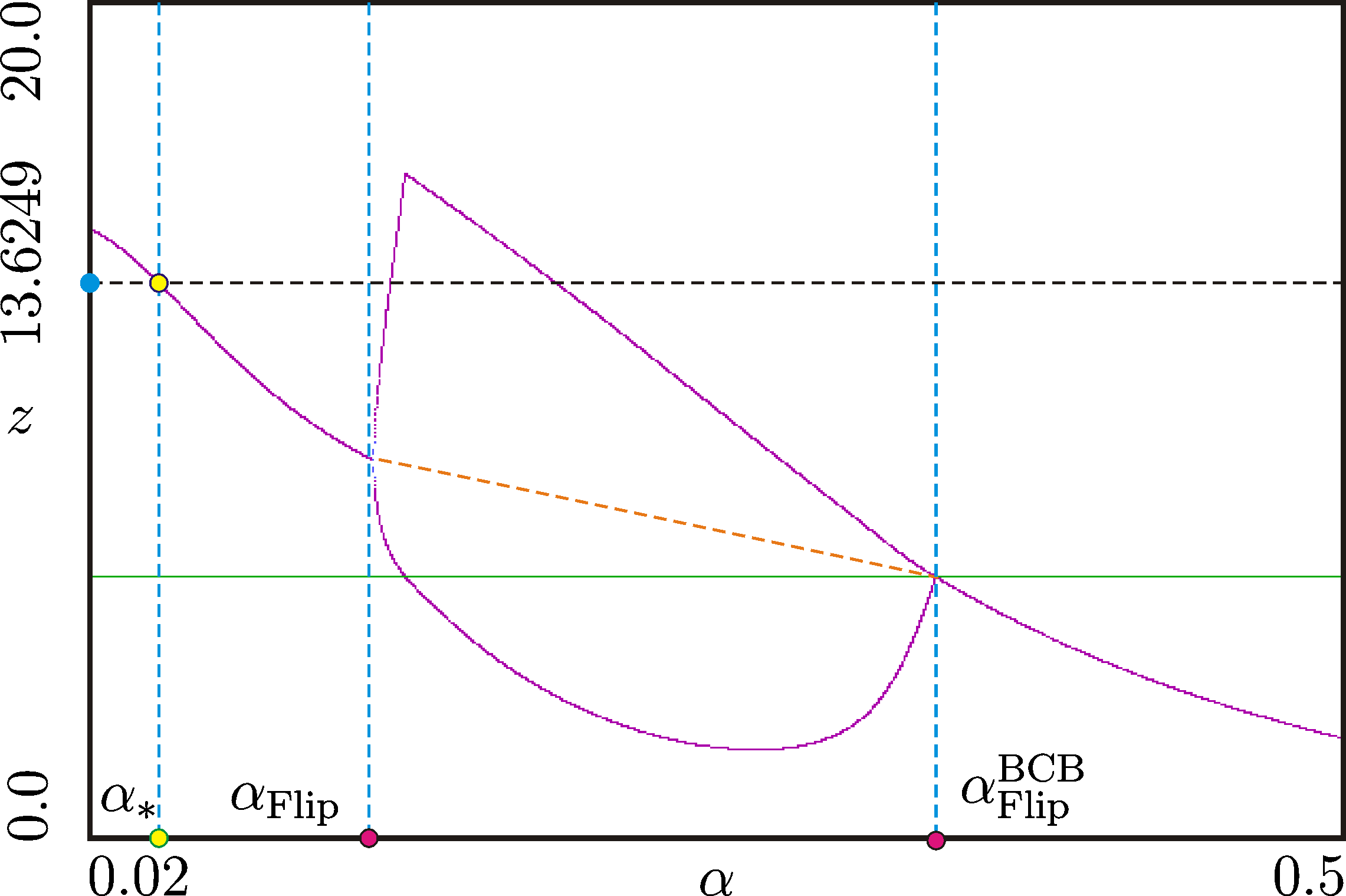}}
\\
\caption{\label{f10}   Bifurcation diagram for $F_1=150.0$, $F_2=400.0$, $\Phi_1=11.0$, $\Phi_2=50.0$. Here  $\varphi(\bar{y}_0)=2.1256$ for $\bar{y}_0=13.6249$ and $\alpha=\alpha_*$, $\alpha_*\approx0.0374\approx \bar\alpha$. The intervals $0.02<\alpha<\alpha_\mathrm{Flip}$ and  $\alpha\gtrapprox\alpha_\mathrm{Flip}^\mathrm{BCB}$ are the regions of the  stability of a fixed point. Here   $\alpha_\mathrm{Flip}\approx 0.107$, $\alpha_\mathrm{Flip}^\mathrm{BCB}\approx 0.277$.
The dotted line marks the unstable fixed point. Between the points $\alpha_\mathrm{Flip}$  and $\alpha_\mathrm{Flip}^\mathrm{BCB}$ exists the stable 2-cycle.}
\end{figure}

A simulation of the designed algorithm performance in the closed-loop for the nominal parameters $(\bar\alpha,\bar\gamma$), starting induction of NMB  ($x(0)=0$, $y(0)=100\%$ ), is depicted in Fig.~\ref{fig:NMB_1-cycle}.  The corridor of the output values under stationary conditions is calculated according to \cite{MPZ24}. The output signal converges promptly to the stationary solution with a minimal overshoot. To demonstrate robustness of the controller with respect to variations in PK, select the plant parameters as $(\alpha=0.11, \bar\gamma$). In Fig.~\ref{fig:NMB_2-cycle}, the transition from ($x(0)=0$, $y(0)=100\%$ ) to stationary solution is shown. Due to the parameter mismatch between the controller and plant, instead of the nominal 1-cycle, a stable 2-cycle arises under stationary conditions, see Fig.~\ref{f10} for bifurcation diagram. Despite an unrealistic increase of $\alpha$ for more than four times of the nominal value, the output is still kept within a reasonable interval $y(t)\in \lbrack 0.27, 12 \rbrack$.

\section{Bifurcation analysis}
This section  presents an overview of  nonlinear dynamics phenomena that can potentially arise in drug dosing automated by a pulse-modulated feedback.

Denote the elements of $X_n=\begin{bmatrix}
    x_{1,n} & x_{2,n} &x_{3,n}
\end{bmatrix}^\intercal$ which gives $T_n= \Phi(x_{3,n})$,    $\lambda_n=
F(x_{3,n})$.  Map \eqref{eq:map} takes the form
\begin{equation}
\begin{bmatrix}\label{eq:map_element}
    x_{1,n+1} \\ x_{2,n+1} \\x_{3,n+1}
\end{bmatrix}
=\e^{AT_n}\begin{bmatrix}
    x_{1,n}+\lambda_n \\ x_{2,n} \\x_{3,n}
\end{bmatrix}.
\end{equation}
An element-wise version of \eqref{eq:map_element} is provided in Appendix.

Then the modulation functions {can be rewritten as}
\begin{gather*}
\Phi(x_{3,n})=
\begin{cases}
\Phi_1,&\;\; 0<x_{3,n}<c_\mathcal{L}^\Phi,\\
\textcolor{black}{\dfrac{100C^\gamma_{50}\:k_2}{C^\gamma_{50}+{x_{3,n}}^\gamma}+k_1},&\;\;
c_\mathcal{L}^\Phi\leqslant
x_{3,n}\leqslant
c_\mathcal{R}^\Phi,\\
{\Phi_2},&\;\;
x_{3,n}>c_\mathcal{R}^\Phi,
\end{cases}
\end{gather*}
\begin{gather*}
\quad F(x_{3,n})=
\begin{cases}
F_1,&\;\; 0<x_{3,n}<c_\mathcal{L}^F,\\
\textcolor{black}{\dfrac{100C^\gamma_{50}\:k_4}{C^\gamma_{50}+{x_{3,n}}^\gamma}+k_3},&\;\;
c_\mathcal{L}^F\leqslant
x_{3,n}\leqslant
c_\mathcal{R}^F,\\
{F_2},&\;\;
x_{3,n}>c_\mathcal{R}^\Phi.
\end{cases}
\end{gather*}

Map \eqref{eq:map_element} is piecewise-smooth due to the discontinuities between the affine segments of the modulation functions. Such maps  are
characterized by the fact that their phase space is divided into
regions with different dynamics, separated from each other
by the so-called switching sets. 
The saturation limits of the modulation functions result in the switching manifolds
$c_\mathcal{L}^\Phi$, $c_\mathcal{R}^\Phi$, $c_\mathcal{L}^F$, $c_\mathcal{R}^F$ 
\begin{gather*}
c_\mathcal{L}^\Phi=c_{50}\cdot\left(\dfrac{100\:k_2}{\Phi_1-k_1}-1\right)^{\dfrac1\gamma},\;
c_\mathcal{R}^\Phi=c_{50}\cdot\left(\dfrac{100\:k_2}{\Phi_2-k_1}-1\right)^{\dfrac1\gamma},\\
c_\mathcal{L}^F=c_{50}\cdot\left(\dfrac{100\:k_4}{F_1-k_3}-1\right)^{\dfrac1\gamma},\;
c_\mathcal{R}^F=c_{50}\cdot\left(\dfrac{100\:k_4}{F_2-k_3}-1\right)^{\dfrac1\gamma}.
\end{gather*}

\begin{figure}[ht]
\centering \centering
\parbox[c]{0.65\linewidth}{\includegraphics[width=\linewidth]{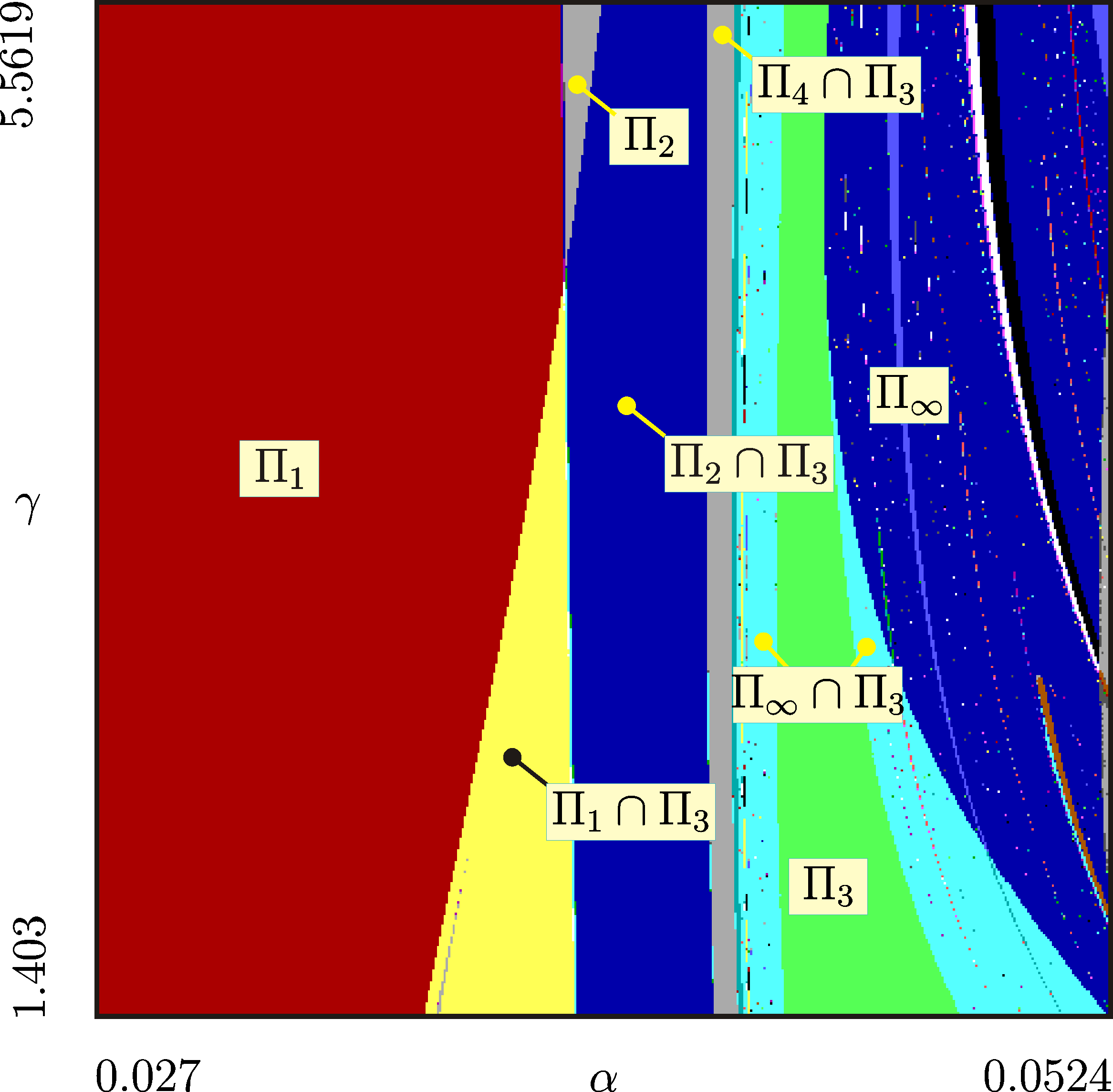}}
\\
\caption{\label{f12} 2D  bifurcation diagram for $F'(y_0)=-1.0$, and $\Phi'(y_0)=4.0$,  $F_1=150.0$, $F_2=400.0$, $\Phi_1=5.0$, $\Phi_2=50.0$, $0.027<\alpha<0.0524$, $1.403<\gamma<5.5619$, $a_1=-0.0374$, $a_2=-0.1496$, $a_3=-0.374$, $\mathrm{g}_1=0.0374$, $\mathrm{g}_2=0.056$, $T=20.0$, $\lambda=300.0$. }
\end{figure}

\begin{figure}[ht]
\centering \centering
\parbox[c]{0.65\linewidth}{\includegraphics[width=\linewidth]{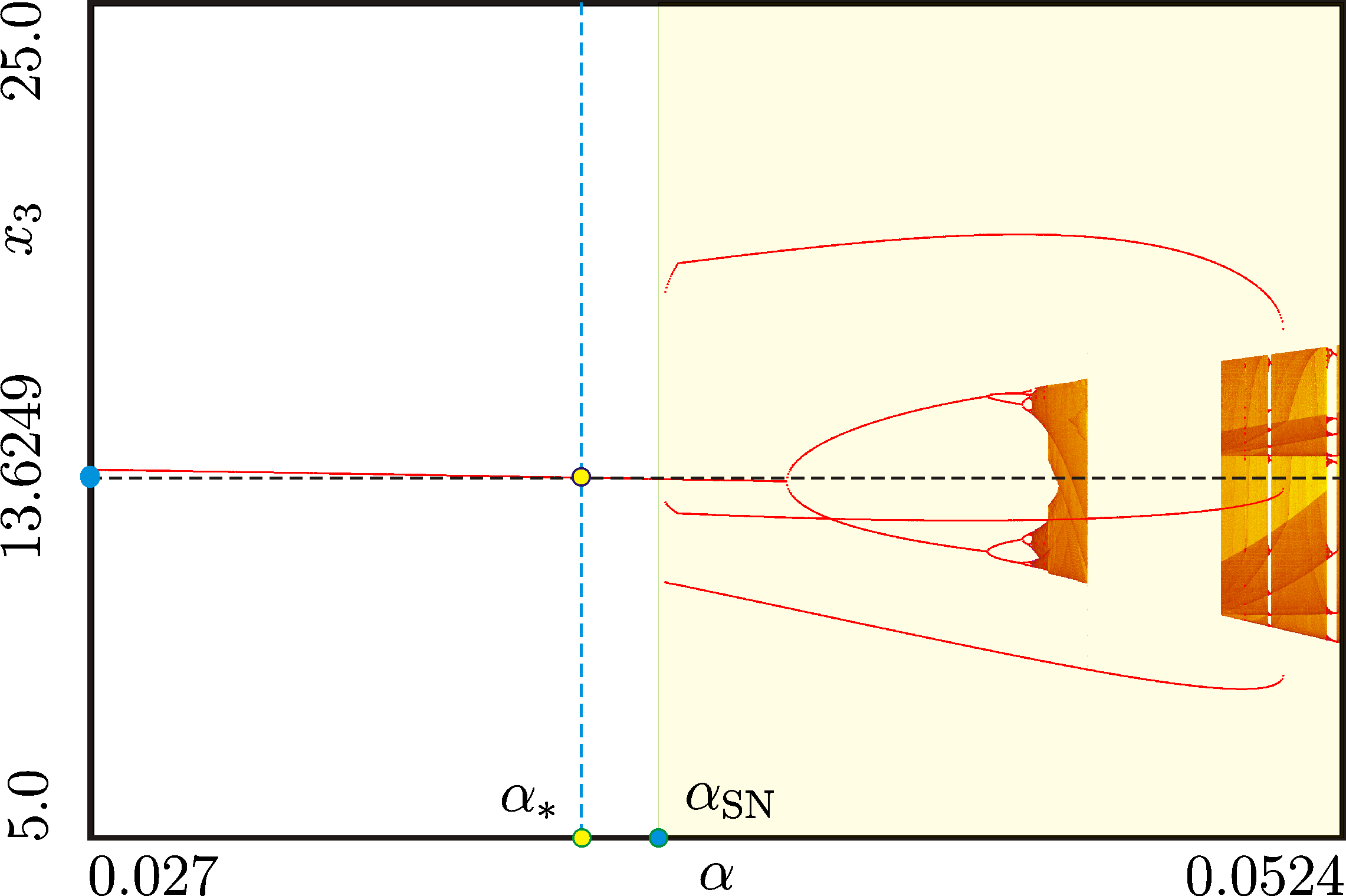}}
\\
\caption{\label{f11} 1D  bifurcation diagram for $F'(y_0)=-1.0$, and $\Phi'(y_0)=4.0$,  $F_1=150.0$, $F_2=400.0$, $\Phi_1=5.0$, $\Phi_2=50.0$, $\gamma=2.6677$, $0.027<\alpha<0.0524$, $a_1=-0.0374$, $a_2=-0.1496$, $a_3=-0.374$, $\mathrm{g}_1=0.0374$, $\mathrm{g}_2=0.056$, $T=20.0$, $\lambda=300.0$,
$k_1=40.87311$, $k_2=-9.819844$, $k_3=294.7817$, $k_4=2.45496$, $\varphi(\bar{y}_0)=2.125605$, $\varphi'(\bar{y}_0)=-0.4073$. 
The multistability region is $\alpha >\alpha_\mathrm{SN} $  (in yellow).
}
\end{figure}

In addition to the bifurcations
occurring in smooth systems, piecewise smooth systems demonstrate a variety of border-collision-related phenomena which
occur when an invariant set such as, for example, a cycle, collides with a switching manifold. This leads to the so-called border-collision bifurcations. An overview
of the phenomena related to border collisions can be found in
~\cite{MZ03}.

The two-parameter
bifurcation diagram in the parameter plane
($\alpha, \gamma$) for $F'(\bar{y}_0)=-1.0$ and $\Phi'(\bar{y}_0)=4.0$ is shown in Fig.~\ref{f12}. Here $\bar{y}_0$ is the third component of the fixed point $\bar{y}_0=13.6249$.  Inspection of this diagram reveals
the presence of a dense set of periodic windows which overlap with other regions of periodicity, for example, $\Pi_1\cap\Pi_3$,  $\Pi_2\cap\Pi_3$, $\Pi_4\cap\Pi_3$, where $\Pi_m$, $m=1,2,3,4$ are the regions of periodicity with period $m$ ($m$-cycles, cf. \eqref{eq:m-cycle}, and $\Pi_\infty$ is the domain of chaotic dynamics.

Depending on the parameter values,  a variety
of different scenarios is observed. 
Fig.~\ref{f11} presents the one-parameter scan for $\gamma=2.6677$. At the point $\alpha_\mathrm{SN}$, a classic  saddle-node  bifurcation occurs, in which the stable and saddle 2-cycles arise. This leads to bistability when a stable fixed point coexists with a stable 3-cycle. 
The basins of attraction for coexisting motions 
are separated by  the stable manifold of the saddle cycle.
As the parameter $\alpha$ increases, the stable 3-cycle undergoes  a  persistence border-collision bifurcation.
With further variation of $\alpha$, one can observe a transition to chaos through a period-doubling sequence (see Fig.~\ref{f11}).
It is well known that numerous windows of periodicity are   found in the region of chaotic dynamics.
Hence  the map displays situations where
several stable cycles (3-cycle and others)  coexist
within a wide range of parameters. These cycles typically arise in   saddle-node or fold border-collision bifurcations
and, with changing parameters, they can undergo sequence of period doubling bifurcations, leading to the transition to chaos.
As a result, there exist parameter domains wherein, alongside with stable cycles,
there are coexisting modes of chaotic oscillations i.e. this leads  to the emergence of multistable dynamics.

A distinguishing feature of multistable systems
is their sensitivity to noise: an arbitrarily small
level of noise may cause a sudden transition from
one attractor to another.

\section{Simulation}
This section illustrates the practical significance  of the complex nonlinear phenomena revealed by bifurcation analysis above by examining time evolutions of the controlled PKPD model output under distinctive simulation scenarios. As in the previous section, let the slopes of the modulation functions be $F'(y_0)=-1.0$, and $\Phi'(y_0)=4.0$, which results in a more aggressive control compared to what is implemented in Section~\ref{sec:design}.  The saturation limits are set to $F_1=150.0$, $F_2=400.0$, $\Phi_1=5.0$, $\Phi_2=50.0$. Then only the value of $\Phi_1$ deviates from that of the nominal design thus allowing for more drug doses to be administered when the output $y(t)$ is far from the fixed-point value $y_0$.
\subsection{Bistability}

Select the PKPD model parameters as $(\alpha=0.04, \bar\gamma)$. According to the bifurcation diagram in Fig.~\ref{f12}, two stable stationary periodic solutions exist, a 1-cycle and a 3-cycle.

\begin{figure}[t]
\centering 
\includegraphics[width=0.9\linewidth]{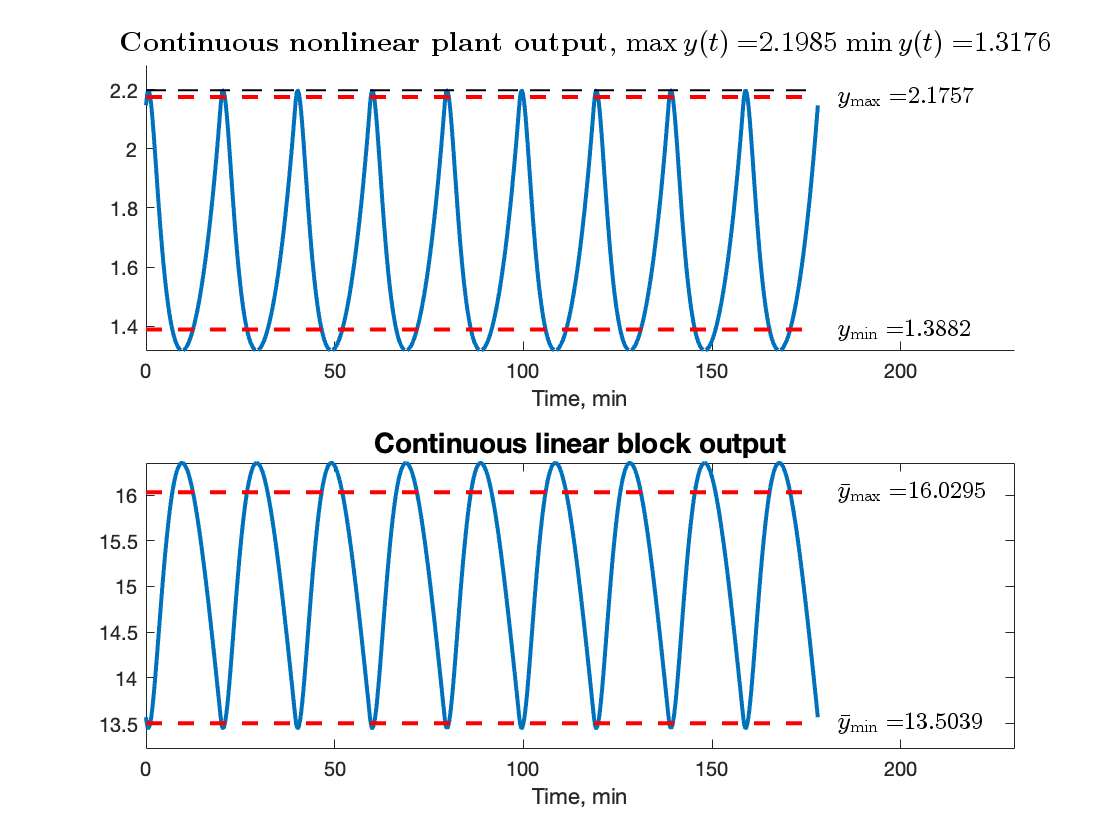}
\caption{ The  1-cycle with $\lambda=300.0549, T=19.7803$ stabilized by the modulation function slopes $F^\prime(y_0)=-1$, $\Phi^\prime(y_0)=4$. The PKPD model parameters are $(\alpha=0.04, \bar\gamma)$ The initial condition on the continuous block is $x(0)=X_{1}$. The desired output corridor values for the output $y(t)$  are  marked by dashed red lines. The actual output corridor is marked by dashed black lines.}\label{fig:1-cycle_bi}
\end{figure}

\paragraph{1-cycle}
 The fixed point of the 1-cycle \[
X_{1}=\begin{bmatrix}
    248.7854486 &78.5195950604 &13.570394339
\end{bmatrix}^\intercal .
\]
The periodic solution corresponding to $X_1$ is shown in Fig.~\ref{fig:1-cycle_bi}. Since $\alpha\ne\bar\alpha$, the actual 1-cycle parameters deviate slightly from the nominal ones and amount to $\lambda=300.0549, T=19.7803$. The changes in the output corridor values are also insignificant. 

\begin{figure}[t]
\centering 
\includegraphics[width=0.9\linewidth]{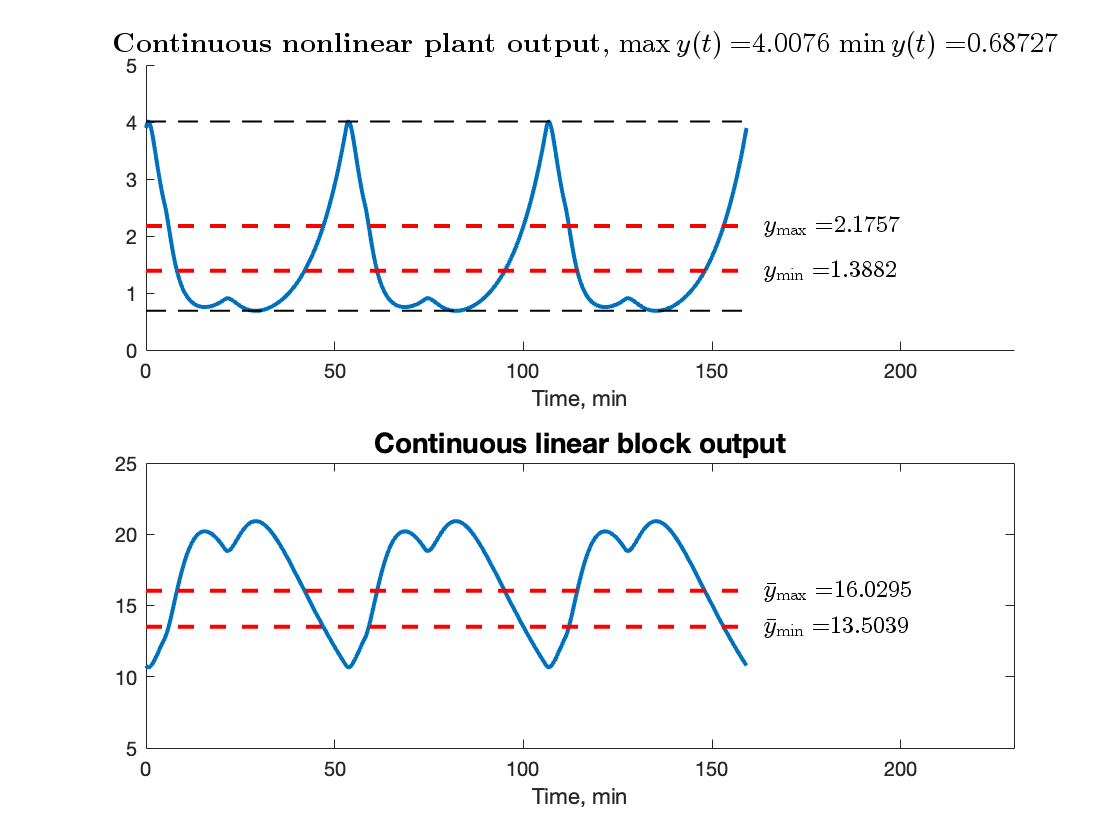}
\caption{ The  3-cycle with $\lambda=300.0549, T=19.7803$ stabilized by the modulation function slopes $F^\prime(y_0)=-1$, $\Phi^\prime(y_0)=4$. The PKPD model parameters are $(\alpha=0.04, \bar\gamma)$ The initial condition on the continuous block is $x(0)=X_{3,1}$. The desired output corridor values for the output $y(t)$  are  marked by dashed red lines. The actual output corridor is marked by dashed black lines.}\label{fig:3-cycle_bi}
\end{figure}

\paragraph{3-cycle}
In Fig.~\ref{fig:3-cycle_bi}, the stable 3-cycle co-existing with the 1-cycle above is depicted. The simulation is initiated from the fixed point
\[
X_{3,1}=\begin{bmatrix}
184.8970537 &60.984174899 &10.784038691
\end{bmatrix}^\intercal.
\]
with the doses and dose times (feedback firings) 
$T_1=5$, $\lambda_1=304.3431$, $T_2=15.9261$, $\lambda_2=301.0185$, $T_3=32.1123$, $\lambda_3=296.9719$. Both the signal form and the actual output corridor are significantly different from what is observed in the 1-cycle. Yet, the overall controller performance is satisfactory also in this operation mode, as the measured output remains within the clinically safe interval but with a tendency to overdose.

\begin{figure}[t]
\centering 
\includegraphics[width=0.9\linewidth]{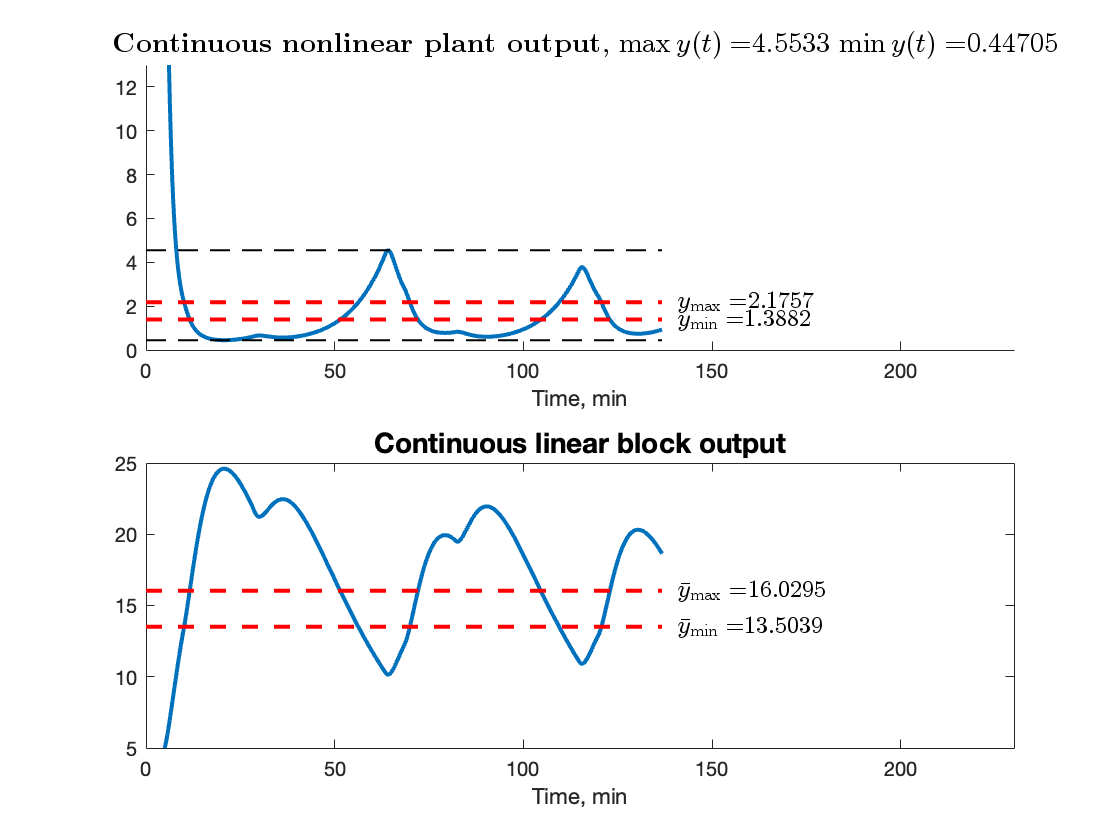}
\caption{ Convergence to the  3-cycle from $x(0)=0$  in the NMB model with $(\alpha=0.04,\bar\gamma)$ stabilized by the modulation function slopes $F^\prime(y_0)=-1$, $\Phi^\prime(y_0)=4$. Top plot: the nonlinear output $y(t)$. The horizontal black dashed lines mark $\inf_t y(t)$ and $\sup_{t\in \lbrack T,5T\rbrack} y(t)$. The stationary output corridor values for  the 1-cycle are marked in red. Bottom plot: the linear output $\bar y(t)$.
}\label{fig:3-cycle_from_zero}
\end{figure}

 In a more realistic simulation scenario, when the administration is initiated from $x(0)=0$, i.e. $y(0)=100\%$, the measured output converges to the 3-cycle as its domain of  attraction is first crossed by the solution, see Fig.~\ref{fig:3-cycle_from_zero}. Therefore, the nominal drug administration regimen (the 1-cycle in Fig.~\ref{fig:1-cycle_bi}) is never reached despite its existence and orbital stability being guaranteed by design. Notice also that the overdosing is worsened in this case compared to the least output value of the output in Fig.~\ref{fig:3-cycle_bi}, from $\inf y(t)=0.68727$ to $\inf y(t)=0.44705$.

\subsection{Deterministic chaos}
As seen in Fig.~\ref{f12}, a further increase in $\alpha$ leads to chaotic  solutions in the closed-loop system. They are characterized by lack of a periodic pattern under stationary conditions. Yet, since all the solutions of closed-loop system \eqref{eq:1_wiener}, \eqref{eq:nonlin_NMB} are bounded (cf. {\bf C3} ), the chaotic solutions are also bounded and stay within an attractor when initiated in it.

\begin{figure}[t]
\centering 
\includegraphics[width=0.9\linewidth]{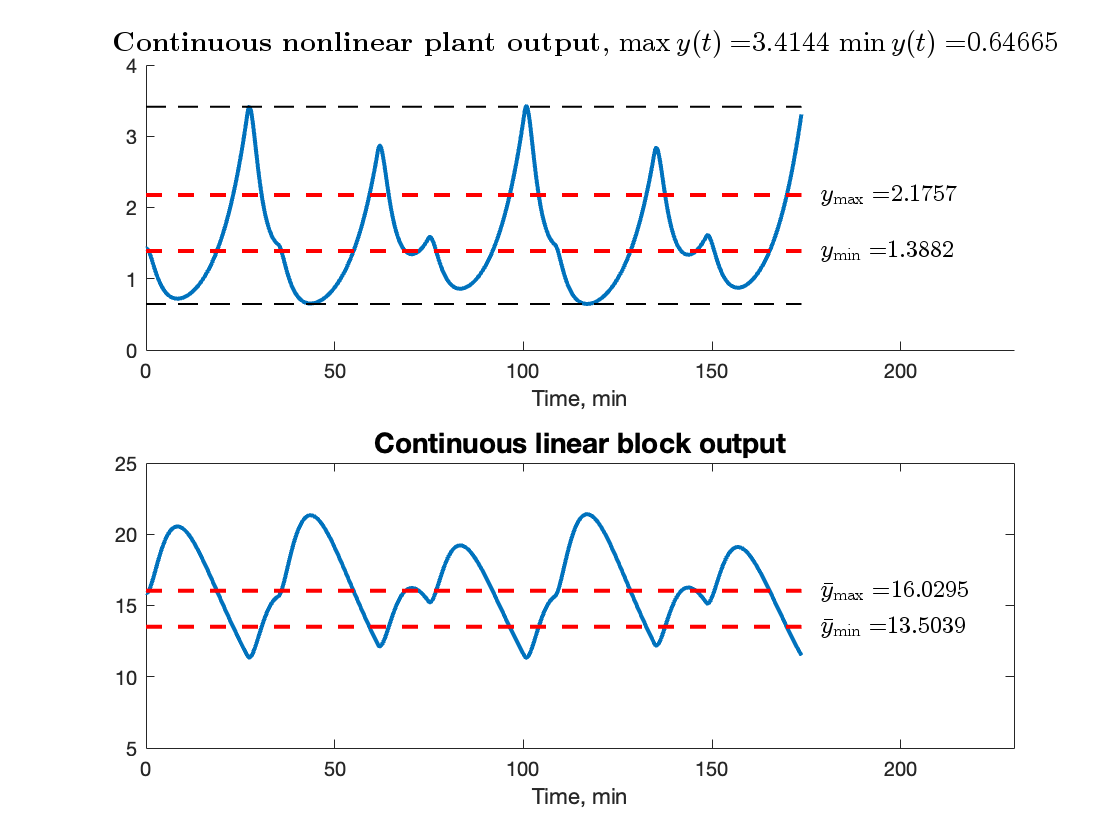}
\caption{ A chaotic solution  in the NMB model with $(\alpha=0.0467,\bar\gamma)$ stabilized by the modulation function slopes $F^\prime(y_0)=-1$, $\Phi^\prime(y_0)=4$. Top plot: the nonlinear output $y(t)$. The horizontal black dashed lines mark $\inf_t y(t)$ and $\sup_{t\in \lbrack T,5T\rbrack} y(t)$. The stationary output corridor values for  the 1-cycle are marked in red. Bottom plot: the linear output $\bar y(t)$.
}\label{fig:chaos}
\end{figure}

Let the PK parameter be $\alpha=0.0467$ and select the initial conditions in the continuous block as $x(0)=X_{\infty}$, where
\[
X_{\infty}=\begin{bmatrix}
295.59510349 &83.573342274 &15.880952316
\end{bmatrix}^\intercal.
\]
A solution that belongs to the chaotic attractor is shown in Fig.~\ref{fig:chaos}. No periodic behavior is observed and the firing parameters of the solution are
$T_1=26.9036$, $\lambda_1=298.2741$, $T_2=8.1826$, $\lambda_2=302.9543$,
$T_3=26.3449$, $\lambda_3=298.4138$, $T_4=13.4515$, $\lambda_4=301.6371$, $T_5=25.4349$, $\lambda_5=298.6413$, $T_6=8.0533$, $\lambda_6=302.9867$, $T_7=26.3354$, $\lambda_7=298.4161$, $T_8=13.7474$, $\lambda_8=301.5632$, $T_9=25.2106$, $\lambda_9=298.6974$.  This particular solution exhibits more variability in the timing of the doses than in their magnitude. Chaotic solutions are highly sensitive to perturbation and should be avoided in engineering practice.

\section*{CONCLUSIONS}
The dynamical properties of a pulse-modulated feedback drug dosing algorithm are studied. Dosing of a neuromuscular blockade agent whose effect is measured by  a train-of-four neuromuscular monitor is selected as illustration. The design of the controller can be performed based on a pharmacokinetic-pharmacodynamic model of the drug and desired dosing regimen. Being designed for population mean parameters, the control law demonstrates high degree of robustness against parameter variability. Since the closed-loop dynamics are highly nonlinear, controller design should be informed by thorough bifurcation analysis in order to safeguard patient safety.




\bibliographystyle{IEEEtran}
\bibliography{observer,refs} 

\begin{thebibliography}{10}
\providecommand{\url}[1]{#1}
\csname url@samestyle\endcsname
\providecommand{\newblock}{\relax}
\providecommand{\bibinfo}[2]{#2}
\providecommand{\BIBentrySTDinterwordspacing}{\spaceskip=0pt\relax}
\providecommand{\BIBentryALTinterwordstretchfactor}{4}
\providecommand{\BIBentryALTinterwordspacing}{\spaceskip=\fontdimen2\font plus
\BIBentryALTinterwordstretchfactor\fontdimen3\font minus
  \fontdimen4\font\relax}
\providecommand{\BIBforeignlanguage}[2]{{%
\expandafter\ifx\csname l@#1\endcsname\relax
\typeout{** WARNING: IEEEtran.bst: No hyphenation pattern has been}%
\typeout{** loaded for the language `#1'. Using the pattern for}%
\typeout{** the default language instead.}%
\else
\language=\csname l@#1\endcsname
\fi
#2}}
\providecommand{\BIBdecl}{\relax}
\BIBdecl

\bibitem{GAH13}
R.~George, T.~Allen, and A.~Habib, ``Intermittent epidural bolus compared with
  continuous epidural infusions for labor analgesia: a systematic review and
  meta-analysis,'' \emph{Anesth Analg.}, vol. 116, pp. 133--144, 2013.

\bibitem{FBW23}
R.~Fernandez~Rojas, N.~Brown, G.~Waddington, and R.~Goecke, ``A systematic
  review of neurophysiological sensing for the assessment of acute pain,''
  \emph{npj Digital Medicine}, vol.~6, pp. 2398--6352, 2023.

\bibitem{RRR22}
J.~Rodríguez-Blanco, T.~Rodríguez-Yanez, J.~D. Rodríguez-Blanco, A.~J.
  Almanza-Hurtado, M.~C. Martínez-Ávila, D.~Borré-Naranjo, M.~C.
  Acuña~Caballero, and C.~Dueñas-Castell, ``Neuromuscular blocking agents in
  the intensive care unit,'' \emph{The Journal of international medical
  research}, vol.~50, no.~9, p. 3000605221128148, 2022.

\bibitem{MH06}
C.~D. McGrath and J.~M. Hunter, ``{Monitoring of neuromuscular block},''
  \emph{Continuing Education in Anaesthesia Critical Care \& Pain}, vol.~6,
  no.~1, pp. 7--12, 02 2006.

\bibitem{TSB21}
M.~L. {Thompson Bastin}, R.~R. Smith, B.~D. Bissell, H.~N. Wolf, A.~M. Wiegand,
  M.~E. Cavagnini, Y.~Ahmad, and A.~H. Flannery, ``Comparison of fixed dose
  versus train-of-four titration of cisatracurium in acute respiratory distress
  syndrome,'' \emph{Journal of Critical Care}, vol.~65, pp. 86--90, 2021.

\bibitem{MCS06}
A.~Medvedev, A.~Churilov, and A.~Shepeljavyi, ``Mathematical models of
  testosterone regulation,'' in \emph{Stochastic optimization in
  informatics}.\hskip 1em plus 0.5em minus 0.4em\relax Saint Petersburg State
  University, 2006, no.~2, pp. 147--158, in Russian.

\bibitem{Aut09}
A.~Churilov, A.~Medvedev, and A.~Shepeljavyi, ``Mathematical model of non-basal
  testosterone regulation in the male by pulse modulated feedback,''
  \emph{Automatica}, vol.~45, no.~1, pp. 78--85, 2009.

\bibitem{SWM12}
M.~M. da~Silva, T.~Wigren, and T.~Mendonca, ``Nonlinear identification of a
  minimal neuromuscular blockade model in anesthesia,'' \emph{IEEE Transactions
  on Control Systems Technology}, vol.~20, no.~1, pp. 181--188, 2012.

\bibitem{ZCM12b}
Z.~T. Zhusubaliyev, A.~Churilov, and A.~Medvedev, ``Bifurcation phenomena in an
  impulsive model of non-basal testosterone regulation,'' \emph{Chaos},
  vol.~22, no.~1, p. 013121, 2012.

\bibitem{PRM24}
A.~V. Proskurnikov, H.~Runvik, and A.~Medvedev, ``Cycles in impulsive
  {G}oodwin’s oscillators of arbitrary order,'' \emph{Automatica}, vol. 159,
  p. 111379, 2024.

\bibitem{MPZ24_ECC}
A.~Medvedev, A.~V. Proskurnikov, and Z.~T. Zhusubaliyev, ``Impulsive feedback
  control for dosing applications,'' in \emph{European Control Conference},
  Stockholm, Sweden, 2024, pp. 1258--1263.

\bibitem{MPZ24}
------, ``Output corridor control via design of impulsive {G}oodwin's
  oscillator,'' in \emph{American Control Conference}, Toronto, Canada, 2024.

\bibitem{MZ03}
E.~Mosekilde and Z.~T. Zhusubaliyev, \emph{Bifurcations and chaos in
  piecewise-smooth dynamical systems}.\hskip 1em plus 0.5em minus 0.4em\relax
  World Scientific, 2003.

\bibitem{DeBoor2005}
C.~De~Boor, ``Divided differences,'' \emph{Surveys in Approximation Theory},
  vol.~1, pp. 46--69, 2005.

\end{thebibliography}
\section*{APPENDIX}


Introduce first divided difference \cite{DeBoor2005} of a function $h(\cdot)$
\[
h\lbrack x_1, x_2 \rbrack \triangleq \frac{h(x_1)-h(x_2)}{x_1-x_2},
\]
and higher-order divided differences defined recursively by 
\[
h\lbrack x_0, \dots, x_k \rbrack = \frac{h\lbrack x_1, \dots, x_k\rbrack-h\lbrack x_0, \dots, x_{k-1}\rbrack }{x_k-x_0}.
\]
The matrix exponential is evaluated in \cite{PRM24} as
\begin{align*}
&\e^{At}=\\
&\left[\begin{array}{ccc}
    \e^{-a_1t} &0 &0\\
    g_1t \e \lbrack -a_1t, -a_2t \rbrack &\e^{-a_2t} &0\\
    g_1g_2t^2 \e \lbrack -a_1t, -a_2t,  -a_3t\rbrack & g_2t\e \lbrack -a_2t, -a_3t \rbrack & \e^{-a_3t}
\end{array}\right], \nonumber
\end{align*}
and \eqref{eq:map_element} takes the form
\begin{align*}
    & x_{1,n+1}= \e^{-a_1T_n}(x_{1,n}+\lambda_n),\\
    & x_{2,n+1}=\e^{-a_2T_n}x_{2,n}+g_1T_n\e \lbrack -a_1t, -a_2t \rbrack(x_{1,n}+\lambda_n),\\
   &  x_{3,n+1}=\e^{-a_3T_n}x_{3,n}+g_2T_n \times \\
     &\big( g_1T_n \e \lbrack -a_1t, -a_2t,  -a_3t\rbrack (x_{1,n}+\lambda_n) + \e \lbrack -a_2t, -a_3t \rbrack x_{2,n} \big).
\end{align*}




\end{document}